\documentclass[12pt,a4paper,leqno,verbatim]{amsart}
\newcounter{minutes}
\divide\time by 60
\newcounter{hours}
\multiply\time by 60 \addtocounter{minutes}{-\time}

\usepackage{amssymb}
\usepackage{hyperref}
\usepackage{graphicx}
\usepackage{subfigure}
\date{}
\newfont{\cyrilic}{wncyr10 scaled 1000}

\title[{Complete $(p,q)$-elliptic integrals}]
{On generalized complete $(p,q)$-elliptic integrals}

\author[B. A. Bhayo]{Barkat Ali Bhayo}
\address{Department of mathematics, Sukkur IBA University, Sindh, Pakistan}
\address{Faculty of Natural Sciences
Sabanci University, 34956 Tuzla/Istanbul, Turkey}
\email{barkat.bhayo@iba-suk.edu.pk}

\author[N. G. G\"o\u{g}\"u\c{s}]{Nihat G\"okhan G\"o\u{g}\"u\c{s}}
\address{Faculty of Natural Sciences
Sabanci University, 34956 Tuzla/Istanbul, Turkey}
\email{nggogus@sabanciuniv.edu}

\author[L. Yin]{Li Yin}
\address{Department of Mathematics, Binzhou University, Binzhou City, Shandong Province, 256603, China}
\email{yinli\_79@163.com}

\newcommand{\comment}[1]{}

\swapnumbers
\theoremstyle{plain}

\newtheorem{theorem}[equation]{Theorem}
\newtheorem{lemma}[equation]{Lemma}

\newtheorem{remark}[equation]{Remark}

\numberwithin{equation}{section}

\begin{document}
\font\fFt=eusm10 
\font\fFa=eusm7  
\font\fFp=eusm5  
\def\K{\mathchoice
{\hbox{\,\fFt K}}
{\hbox{\,\fFt K}}
{\hbox{\,\fFa K}}
{\hbox{\,\fFp K}}}
\def\E{\mathchoice
{\hbox{\,\fFt E}}
{\hbox{\,\fFt E}}
{\hbox{\,\fFa E}}
{\hbox{\,\fFp E}}}

\allowdisplaybreaks

\begin{abstract}
In this paper authors study the generalized complete $(p,q)$-elliptic integrals of the first and the second kind as an application of generalized trigonometric functions with two parameters, and establish the Tur\'an type inequalities of these functions.
\end{abstract}

\vspace{.5cm}

\def\thefootnote{}
\footnotetext{ \texttt{\tiny File:~\jobname .tex,
          printed: \number\year-0\number\month-\number\day,
          \thehours.\ifnum\theminutes<10{0}\fi\theminutes}
} \makeatletter\def\thefootnote{\@arabic\c@footnote}\makeatother

\maketitle



\section{Introduction}
For given complex numbers $a,b$ and $c$ with $c\neq0,-1,-2,\ldots$,
the \emph{Gaussian hypergeometric function} ${}_2F_1$ is the
analytic continuation to the slit place $\mathbb{C}\setminus[1,\infty)$ of the series
   $$F\left(a,b;c;z\right) = {}_2F_1\left(a,b;c;z\right)=\sum_{n\geq 0}\frac{(a)_n(b)_n}
                             {(c)_n}\frac{z^n}{n!},\qquad |z|<1.$$
Here $(a)_n$ is the Pochhammer symbol (rising factorial) $(\cdot)_n:\mathbb{C}\to \mathbb{C}$, defined by
   $$(z)_n = \frac{\Gamma(z+n)}{\Gamma(z)} = \prod_{i=1}^n(z+i-1)$$
for $n\in\mathbb{Z}$, see \cite{sl}. 

The integral representation of the hypergeometric function is given as follows \cite[p. 20]{sl}
\begin{equation}\label{hyp}
F(a,b;c;z)=\frac{\Gamma(c)}{\Gamma(b)(c-b)}\int_0^1{t^{b-1}(1-t)^{c-b-1}(1-zt)^{-a}dt}
\end{equation}
${\rm Re}(c)>{\rm Re}(b)>0, |\arg(1-z)|<\pi$.

Special functions, such as the classical \emph{gamma function} $\Gamma$, the {\it digamma function} $\psi$ and the \emph{beta function}  $B(\cdot,\cdot)$ have close relation with hypergeometric function. For $x,y>0$, these functions are defined by
   $$\Gamma(x) = \int^\infty_0 e^{-t}t^{x-1}\,dt,\quad
       \psi(x) = \frac{\Gamma'(x)}{\Gamma(x)},\quad
        B(x,y) = \frac{\Gamma(x)\Gamma(y)}{\Gamma(x+y)},$$
respectively.

For $x \in (0, 1)$ and $p > 1$ the generalized inverse trigonometric functions are defined as follows
$$\arcsin_p(x)=\int_0^x(1-t^p)^{-1/p}dt=x\,F \left(\frac{1}{p},\frac{1}{p};1+\frac{1}{p};x^p\right),$$
$$\arctan_p(x)=\int_0^x(1+t^p)^{-1}dt=x\,F \left(1,\frac{1}{p};1+\frac{1}{p};-x^p\right),$$
and $\arccos_p(x)=\arcsin_p((1-x^p)^{1/p})$, see \cite{bbt,bbv,bvarxiv}. We note that the eigenvalue problem \cite{dm}, $1 < p < \infty$
$$-\Delta_pu=-(|u'|^{p-2}u')'=\lambda |u|^{p-2}u,\quad u(0)=u(1)=0,$$
has eigenvalues $\lambda_n=(p-1)(n\pi_p)^p$ and eigenfunctions $t\mapsto \sin_p(n\pi_p t),\,n\in\mathbb{N}$, where $\sin_p$ is the inverse
function of $\arcsin_p$ and 
$$\pi_p=\frac{2}{p}\int_0^1(1-s)^{-1/p}s^{1/p-1}ds=\frac{2}{p}B\left(1-\frac{1}{p},\frac{1}{p}\right)=\frac{2\pi}{p\,\sin(\pi/p)}.$$ 
The other generalized trigonometric functions 
$\cos_p : (0; \pi_p/2) \to (0, 1)$ and $\tan_p : (0, b_p) \to
(0, 1)$ are defined as the inverse of the generalized inverse
trigonometric functions $\arccos_p$ and $\arctan_p$, respectively, where
$$b_p=2^{-1/p}\,F \left(\frac{1}{p},\frac{1}{p};1+\frac{1}{p};\frac{1}{2}\right).$$

Recently, Takeuchi \cite{t} studied the $(p,q)$-trigonometric functions depending on two parameters. For $p=q$, these functions 
reduce to the so-called $p$- trigonometric functions introduced by Lindqvist in his highly cited paper \cite{l}. In present, there has been a vivid interest on the generalized trigonometric and hyperbolic functions, numerous papers have been published on the studies of generalized trigonometric functions and their inequalities, see, e.g., \cite{bbv, bsand, bvjapprox, bvarxiv,egl} and the references therein. 

The following ($p,q$)-eigenvalue problem with Dirichl\'et boundary condition was considered by Dr\'abek and Man\'asevich 
\cite{dm}.
Let $\phi_p(x)=|x|^{p-2}x.$ For $T,\lambda>0$ and $p,q>1$
$$\left\{\begin{array}{lll}\displaystyle(\phi_p(u'))'+
\lambda\,\phi_q(u)=0,\quad t\in(0,T),\\
                   \displaystyle u(0)=u(T)=0.\end{array}\right.$$
They 
found the complete solution to this problem. The solution of this problem also 
appears in \cite[Thm 2.1]{t}.
In particular, for $T=\pi_{p,q}$ the function $u(t)=\sin_{p,q}(t)$ is
a solution to this problem with $\lambda=q(p-1)/p$, where
$$\pi_{p,q}=2\int_0^1{(1-t^q)^{-1/p}}\,dt=
\frac{2}{q}B\left(1-\frac{1}{p},\frac{1}{q}\right).$$
For $p=q$, $\pi_{p,q}$ reduces to $\pi_p$, see, e.g., \cite{bbv}.
In order to give the definition of the
function $\sin_{p,q}\,,$ first we define its inverse function ${\rm arcsin}_{p,q}\,,$ then the function itself.
For $x\in[0,1]$, set
$$F_{p,q}(x)={\rm arcsin}_{p,q}(x)=\int_0^x{(1-t^q)}^{-1/p}\,dt\,.$$
The function $F_{p,q}:[0,1]\to [0,\pi_{p,q}/2]$ is an increasing homeomorphism, and 
$$\sin_{p,q}= F_{p,q}^{-1}\,$$
is defined on the the interval $[0,\pi_{p,q}/2]$. The function
$\sin_{p,q}$ can be extended to $[0,\pi_{p,q}]$ by 
$$\sin_{p,q}(x)=\sin_{p,q}(\pi_{p,q}-x),
\quad x\in[\pi_{p,q}/2,\pi_{p,q}].$$
By oddness, the further extension can be made to $[-\pi_{p},\pi_{p}]$. Finally the functions $\sin_{p,q}$ is extended to whole $\mathbb{R}$ by $2\pi_{p}$-periodicity, see \cite{egl}.

The generalized cosine function $\cos_{p,q}$ can be defined as
$$\cos_{p,q}(x)=\frac{d}{dx}\sin_{p,q}(x),\quad x\in \mathbb{R}.$$
One can see easily that $\cos_{p,q}$ is even with period $2\pi_{p,q}$ and odd in about $\pi_{p,q}/2$. Setting $y=\sin_{p,q}(x)$ and letting $x\in[0,\pi_{p,q}/2]$, we get
$$\cos_{p,q}(x)=(1-y^q)^{1/p}=(1-\sin_{p,q}(x)^q)^{1/p}.$$
Clearly, $\cos_{p,q}$ is strictly decreasing with 
$\cos_{p,q}(0)=1$ and $\cos_{p,q}(\pi_{p,q}/2)=0$. From the above definition it follows that
\begin{equation}\label{sincos}
|\cos_{p,q}(x)|^p+|\sin_{p,q}(x)|^q=1,\quad x\in\mathbb{R}.
\end{equation}
The generalized tangent function $\tan_{p,q}$ is defined by
$$\tan_{p,q}(x)=\frac{\sin_{p,q}(x)}{\cos_{p,q}(x)},\quad x\in\{\mathbb{R}: x\neq (z+1/2)\pi_{p,q},\,z\in\mathbb{Z}\}.$$
The usual elementary trigonometric functions are the special case of these $(p,q)$-trigonometric functions when $p=q=2$.

Before we state our main result, we define the generalized complete $(p,q)$-elliptic integrals 
motivated by the work of Takeuchi \cite{t,t1,t2}.
For $p,q>1$, $x\in[0,1]$ and $r\in[0,1)$, we define the generalized complete $(p,q)$-elliptic integral of the first kind $K_{p,q}$ as follows:
$$K_{p,q}(r)=\displaystyle\int_0^1{\frac{dt}{\sqrt[p]{(1-t^q)(1-r^qt^q)^{p-1}}}}=
\displaystyle\int_0^{\pi_{p,q}/2}\frac{dt}{{(1-r^q\sin_{p,q}(t)^q)^{1-1/p}}}.$$
For $p=q$, the function $K_{p,q}$ coincides with $K_p$, where $K_p$ is generalized complete $p$-elliptic integral of the first kind, 
see \cite{t,t1}.

For $p,q>1$, $x\in[0,1]$ and $r\in[0,1)$, we define the function ${\rm arcsn}_{p,q}:[0,1]\to[0,\hat{K}_{p,q}]$ by
$${\rm arcsn}_{p,q}(x)={\rm arcsn}_{p,q}(x,r)=\int_0^x{\frac{dt}{\sqrt[p]{(1-t^q)(1-r^qt^q)}}},$$
and call it generalized inverse Jacobian elliptic function, where $\hat{K}_{p,q}$ is also called generalized complete $(p,q)$-elliptic integral of the first kind, and defined as
\begin{equation}\label{kpq}
\hat{K}_{p,q}(r)=\int_0^1{\frac{dt}{\sqrt[p]{(1-t^q)(1-r^qt^q)}}}={\rm arcsn}_{p,q}(1,r).
\end{equation}
The integral
$\hat{K}_{p,q}$ appears in the study of bifurcation problems of $p$-Laplacian, see \cite{t}.
Substituting $t=\sin_{p,q}(\theta)$ in formula \eqref{kpq}, we get
$$\hat{K}_{p,q}(r)=\int_0^{\pi_{p,q}/2}\frac{d\theta}{\sqrt[p]{1-r^q\sin_{p,q}(\theta)^q}}.$$
Clearly, ${\rm arcsn}_{p,q}(x,r)$ is strictly increasing in $x$, and its inverse
${\rm sn}_{p,q}:[0,\hat{K}_{p,q}]\to[0,1]$ is also strictly increasing, and called generalized Jacobian elliptic function \cite{t}.

Letting ${t=x^{1/q}}$ in \eqref{kpq} and utilizing the formula \eqref{hyp}, we get
\begin{eqnarray*}
\hat{K}_{p,q}(r)&=&\frac{1}{q}\int_{0}^1{x^{1/q-1}(1-x)^{-1/p}(1-xr^q)^{-1/p}dx}\\
&=&\frac{1}{q}\frac{B(1/q,1-1/p+1/q-1/q)}{B(1/q,1-1/p)}\\
& &\times\int_{0}^1{x^{1/q-1}(1-x)^{1-1/p+1/q-1/q-1}(1-xr^q)^{-1/p}dx}\\
&=& \frac{1}{q}B\left(\frac{1}{q},1-\frac{1}{p}\right)
F\left(\frac{1}{p},\frac{1}{q};1-\frac{1}{p}+\frac{1}{q};r^q\right)\\
&=&\frac{\pi_{p,q}}{2}F\left(\frac{1}{p},\frac{1}{q};1-\frac{1}{p}+\frac{1}{q};r^q\right).
\end{eqnarray*}

We define the generalized complete $(p,q)$-elliptic integral of the second kind by
\begin{equation}\label{Ellpq}
E_{p,q}(r)=\int_0^{\pi_{p,q}/2}\sqrt[p]{1-r^q
\sin_{p,q}(t)^q}dt.
\end{equation}
Substituting $x=\sin_{p,q}(t)$ and applying \eqref{sincos}, 
the formula \eqref{Ellpq} can be written as
\begin{equation}\label{epq}
E_{p,q}(r)=\int_{0}^1{(1-x^q)^{-1/p}(1-x^qr^q)^{1/p}dx}.
\end{equation}
By applying formulas \eqref{epq} and \eqref{hyp}, the function $E_{p,q}(r)$ can be expressed in terms of hypergeometric function as below 
$$E_{p,q}(r)=\frac{\pi_{p,q}}{2}F\left(-\frac{1}{p},\frac{1}{q};1-\frac{1}{p}+\frac{1}{q};r^q\right).
$$

\begin{figure}[h]
\includegraphics[width=12cm]{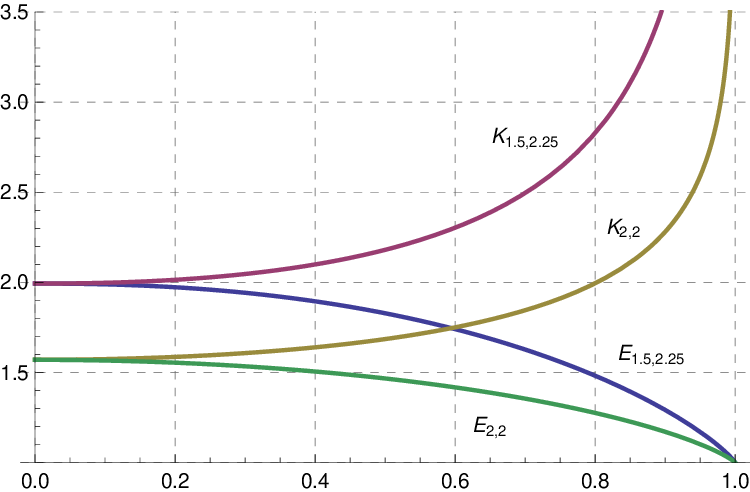}
\caption{Graphs of $K_{1.5,2.25}(r),\,E_{1.5,2.25}(r),\,K_{2,2}(r)=\K(r)$ and $E_{2,2}(r)=\E(r)$, with 
$K_{1.5,2.25}(1)=E_{1.5,2.25}(1) =\pi_{1.5,2.25}/2 \approx 1.9937$.}
\end{figure}

For convenience of the reader we recall the definition of generalized complete $(p,q)$-elliptic integrals of the first and the second kind
as bellow, for $p,q>1$ and $r\in(0,1)$, 

\begin{equation}\label{KEpq}
\left\{\begin{array}{lll} K_{p,q}(r)=\displaystyle\int_0^{\pi_{p,q}/2}\frac{dt}{(1-r^q\sin_{p,q}(t)^q)^{1-1/p}}=
\displaystyle\int_0^1{\frac{dt}{\sqrt[p]{(1-t^q)(1-r^qt^q)^{p-1}}}}\\

\bigskip

 E_{p,q}(r)=\displaystyle\int_0^{\pi_{p,q}/2}\sqrt[p]{1-r^q
\sin_{p,q}(t)^q}dt=\displaystyle\int_0^1\sqrt[p]{\frac{1-r^qt^q}{1-t^q}}dt\\
                      K_{p,q}(0)=\displaystyle\frac{\pi_{p,q}}{2}=
											E_{p,q}(0),\, K_{p,q}(1)=\infty,
											\,E_{p,q}(1)=1.
\end{array}\right.
\end{equation}

For $p=q$, we write $K_{p,p}(r)=K_{p}(r)$ and $E_{p,p}(r)=E_{p}(r)$. Applying formula \eqref{hyp} to \eqref{KEpq}, we get easily
$$K_{p}(r)= \frac{\pi_{p}}{2}F\left(1-\frac{1}{p},\frac{1}{p};1;r^p\right),$$
and 
$$E_{p}(r)= \frac{\pi_{p}}{2}F\left(-\frac{1}{p},\frac{1}{p};1;r^p\right).$$

 Obviously $K_p=\K$ and $\,E_p=\E$ for $p=2$, where $\K$ and $\E$
are the classical elliptic integrals of the first and second kind, respectively. We refer to reader to see \cite{avvb} for the history and the large study on these functions. The generalization and the inequalities of elliptic integrals were studied by numerous authors after the publication of the landmark paper \cite{bbg}. 
The generalized elliptic integrals of the first and second kind on $(0,1)$ are defined respectively by
$$\K_{a,b,c}(r)=\frac{B(a,b)}{2}F(a,b;c;r^2),$$
$$\E_{a,b,c}(r)=\frac{B(a,b)}{2}F(a-1,b;c;r^2),$$
for $0<a<\min\{c,1\}$ and $0<b<c\leq a+b$, see \cite{aqvv}. Clearly,
$\K_{1/2,1/2,1}=\K,\,\E_{1/2,1/2,1}=\E$. For the monotonic properties and the inequalities of these functions, see \cite{aqvv,hlvv}. 
For the historical background and study about the classical and generalized case of elliptic integrals, we refer to reader to see, e.g., 
\cite{aqvv, avvb, b2204, bv2204, hlvv}, and the bibliography of these papers.

Before we present the main results of this paper we recall some definitions as follows:

 A function $f \colon (0,\infty)\to(0,\infty)$ is said to be logarithmically convex, or $\log$-convex, if its natural logarithm $\ln f$ is convex, that is, for all $x,y>0$ and $\lambda\in[0,1]$ we have
   $$f(\lambda x+(1-\lambda)y) \leq \left[f(x)\right]^{\lambda}\left[f(y)\right]^{1-\lambda}.$$
The function $f$ is $\log$-concave if the above inequality is reversed.  

A function $g \colon (0,\infty)\rightarrow(0,\infty)$ is said to be geometrically (or multiplicatively)
convex if it is convex with respect to the geometric mean, that is, if for all $x,y>0$ and all $\lambda\in[0,1]$ the inequality
   $$g(x^{\lambda}y^{1-\lambda}) \leq[g(x)]^{\lambda}[g(y)]^{1-\lambda}$$
holds. The function $g$ is called geometrically concave if the above inequality is reversed. We also note that the differentiable function $f$ is $\log$-convex ($\log$-concave) if and only if
$x \mapsto f'(x)/f(x)$ is increasing (decreasing), while the differentiable function $g$ is geometrically convex (concave) if
and only if the function $x \mapsto xg'(x)/g(x)$ is increasing (decreasing), for more details see
\cite{baricz}. If functions $f$ and $g$ are related by $f(x)=g(e^x)$, then the logarithmic convexity (concavity) of $f$ is equivalent to the geometric convexity (concavity) of $g$.

This paper consists of three sections. In the first section, we give the introduction and definitions of the functions which are being studied in the paper. In the second section, we state our main results. The third section contains few lemmas and the proof of the main result.
	
	\section{Main result}
	
	The main result of this paper reads as follows.
	
	\begin{theorem}\label{thm1}
For $p,q>1,\,r\in(0,1)$ and 
$r'   = (1 - r^q)^{1/q}$, we have
\begin{equation}\label{Kd}
\frac{d}{{dr}}{K_{p,q} (r)}  = \frac{1}{{r\left( {r'} \right)^q }}\left( {E_{p,q} (r) - \left( {r'} \right)^q K_{p,q} (r)} \right)
\end{equation}

\begin{equation}\label{Ed}
\frac{d}{{dr}}\left( {E_{p,q} (r)} \right) = \frac{q}{{pr}}\left( {E_{p,q} (r) - K_{p,q} (r)} \right)
\end{equation}

\begin{equation}\label{Kdd}
\frac{{d^2 }}{{dr^2 }} {K_{p,q} (r)}  = \left( {\frac{{q\left( {r'} \right)^q }}{p} + qr^q  - 2\left( {r'} \right)^q } \right)E_{p,q} (r) + \left( {2\left( {r'} \right)^{2q}  - \frac{q}{p}\left( {r'} \right)^q } \right)K_{p,q} (r)
\end{equation}

\begin{equation}\label{Edd}
\frac{{d^2 }}{{dr^2 }} {E_{p,q} (r)}  = \frac{q}{{pr^2 }}\left( {\left( {\frac{q}{p} - \frac{1}{{\left( {r'} \right)^q }} - 1} \right)E_{p,q} 
(r) + \left( {2 - \frac{q}{p}} \right)K_{p,q} (r)} \right).
\end{equation}
\end{theorem}

It is worth to mention that the formula \eqref{Kd} and \eqref{Ed} recently appeared in \cite[Proposition 2.1]{kamiya}. 

\begin{theorem}\label{thmturan} For $p,q>1$ and $r\in(0,1)$, we have
\begin{enumerate}
\item the function $r\mapsto K_{p,q}(r)$ and $r\mapsto \hat{K}_{p,q}(r)$ are strictly increasing and log-convex. Moreover, $r\mapsto K_{p,q}(r)$ and $r\mapsto \hat{K}_{p,q}(r)$ are strictly geometrically convex on $(0,1)$.
\item The function $r\mapsto E_{p,q}(r)$ is strictly decreasing and geometrically concave on $(0,1)$.
\end{enumerate}
\end{theorem}

\begin{theorem}\label{thm2806} For fixed $r\in(0,1)$ and $p>0$,
\begin{enumerate}
\item the function $q\mapsto \hat{K}_{p,q}(r)$ is strictly decreasing and $\log$-convex on $(0,\infty)$,
\item the function $q\mapsto E_{p,q}(r)$ is strictly decreasing and $\log$-convex on $(0,\infty)$.
\end{enumerate}
In particular, for $r\in(0,1)$, the following Tur\'an type inequalities hold true
$$\hat{K}_{p,q}(r)^2\leq \hat{K}_{p,q-1}(r)\hat{K}_{p,q+1}(r),\qquad p>0,\,q>1,$$
$$E_{p,q}(r)^2\leq E_{p,q-1}(r)E_{p,q+1}(r),\qquad p>0,\,q>1.$$
\end{theorem}

\begin{remark} In the above theorem one can establish the counterpart inequalities for the function $K_{p,q}(r)$ that authors were unable to 
prove. \end{remark}

\begin{theorem}\label{yinthm} For $p,q>1$ and $r\in(0,1), \lambda<\frac{1}{2}$, we have

\begin{equation}\label{Kse}
K_{p,q} (r) = \frac{{\pi _{p,q} }}{2}\sum\limits_{n = 0}^\infty
 {\frac{1}{p} - 1 \choose n}\frac{1}{{(1 - \lambda )^{n + 1 - \frac{1}{p}} }} \sum\limits_{j = 0}^n (-1)^{j}{n \choose j}{- \frac{1}{q}\choose j}{
 \frac{1}{p} - 1 - \frac{1}{q} 
 \choose j}^{-1}\lambda ^{n - j} r^{qj} .
\end{equation}

\begin{equation}\label{Khat-se}
\hat{K}_{p,q} (r) = \frac{{\pi _{p,q} }}{2}\sum\limits_{n = 0}^\infty
 {-\frac{1}{p} \choose n}\frac{1}{{(1 - \lambda )^{n +  \frac{1}{p}} }} \sum\limits_{j = 0}^n (-1)^{j}{n \choose j}{- \frac{1}{q}\choose j}{
 \frac{1}{p} - 1 - \frac{1}{q} 
 \choose j}^{-1}\lambda ^{n - j} r^{qj} .
\end{equation}

\begin{equation}\label{eqse}
E_{p,q} (r) = \frac{\pi _{p,q} }{2}\sum\limits_{n = 0}^\infty {
 \frac{1}{p} \choose n}\frac{1}{{(1 - \lambda )^{n - \frac{1}{p}} }} \sum\limits_{j = 0}^n {( - 1)^j } {n \choose j}{ - \frac{1}{q} \choose j}{\frac{1}{p} - 1 - \frac{1}{q} \choose
 j}^{-1}\lambda^{n - j} r^{qj} .
\end{equation}
\end{theorem}

\begin{remark} \rm
If we let $\lambda=0$, then \eqref{Kse} becomes
\begin{eqnarray*}
 K_{p,q} (r) &=& \frac{\pi _{p,q} }{2}\sum\limits_{n = 0}^\infty  
 {\frac{1}{p} - 1 \choose n}
{n \choose n}
( - 1)^n {- \frac{1}{q} \choose n}
 {\frac{1}{p} - 1 - \frac{1}{q} \choose n}^{-1}
r^{qn}   \\
  &=& \frac{{\pi _{p,q} }}{2}\sum\limits_{n = 0}^\infty  {\frac{{\left( {1 - \frac{1}{p}} \right)_n \left( {\frac{1}{q}} \right)_n (1)_n }}
	{\left( {1 + \frac{1}{q} - \frac{1}{p}} \right)_n{ (1)_n }}\frac{{r^{qn} }}{{n!}}}  \\
 & =& \frac{{\pi _{p,q} }}{2}{}_2F_1 \left( {1 - \frac{1}{p},\frac{1}{q};1 + \frac{1}{q} - \frac{1}{p};r^q } \right).
 \end{eqnarray*}

Similarly, by letting $\lambda =0$ in \eqref{Khat-se} and \eqref{eqse} we get
$$
 \hat{K}_{p,q} (r) = \frac{{\pi _{p,q} }}{2}{}_2F_1 \left( {\frac{1}{p},\frac{1}{q};1 + \frac{1}{q} - \frac{1}{p};r^q } \right),
 $$
and 
$$E_{p,q} (r)= \frac{{\pi _{p,q} }}{2}{}_2F_1 \left( { - \frac{1}{p},\frac{1}{q};1 + \frac{1}{q} - \frac{1}{p};r^q } \right),$$
respectively.
\end{remark}

\section{Preliminaries and proofs}
Before we quote here few lemmas which will be used in the proof of theorems.
	
The following lemma follows easily from the definition
	and \eqref{sincos}.
\begin{lemma}\cite[Proposition 3.1]{egl}\label{lemma1}
For all $x\in[0,\pi_{p,q}/2]$,
\begin{enumerate}
\item $\displaystyle\frac{d}{dx}\sin_{p,q}x=\cos_{p,q}x,$\\
\item $\displaystyle\frac{d}{{dx}}\cos _{p,q} x =  - \frac{q}{p}\left( {\cos _{p,q} x} \right)^{2 - p} \left( {\sin _{p,q} x} \right)^{q - 1} ,$\\
\item $\displaystyle\frac{d}{{dx}}\left( -({  \cos _{p,q} x} \right)^{p - 1})  = \frac{{(p - 1)q}}{p}\left( {\sin _{p,q} x} \right)^{q - 1} ,$
\item $\displaystyle\frac{d}{{dx}}\left( {\sin _{p,q} x} \right)^q  = q\left( {\sin _{p,q} x} \right)^{q - 1} \cos _{p,q} x.$
\end{enumerate}
\end{lemma}

\begin{lemma}\label{lemma2}\cite[Lemma 2]{bbv}
If the function $v\mapsto K(v,t)$ is positive and (strictly) geometrically convex on $[a,b]$ for $t\in(0,x)$, with $0<a<b$ and $x>0$. Then the function
$$v\mapsto f_v(x)=\int_0^x{K(v,t)}dt$$
is also (strictly) geometrically convex.
\end{lemma}

\begin{lemma}\cite{ar}
Suppose that the function $G$ given by 
$G(x) = \frac{{g(x)}}{{\left( {1 - \alpha x^\eta  } \right)^\xi  }}$ satisfies $g,G\in L^{1} [0,1]$
where $0<\alpha \leq1, \eta>0, 
\xi\in\mathbb{R}$, and $$b_j  = b_j (\alpha ,\eta ) = \alpha ^j \int_0^1 {t^{j\eta } g(t)dt,} $$
it follows that 
\begin{equation}\label{areq}
\int_0^1 {\frac{{g(x)}}{{\left( {1 - \alpha x^\eta  } \right)^\xi  }}dx = \sum\limits_{n = 0}^\infty  {\frac{{\left( \xi  \right)_n }}{{n!(1 - \lambda )^{n + \xi } }}} } \sum\limits_{j = 0}^n {n\choose j} ( - \lambda )^{n - j} b_j (\alpha,\eta).
\end{equation}
\end{lemma}


\noindent{\bf Proof of Theorem \ref{thm1}.}
Applying the derivative formulas given in Lemma \ref{lemma1} and utilizing the identity \eqref{sincos}, we get 
\begin{equation}\label{equation2004a}
\frac{d}{{dx}}\left( {\frac{{ - \left( {\cos _{p,q} x} \right)^{p - 1} }}{{\left( {1 - k^q \left( {\sin _{p,q} x} \right)^q } \right)^{1 - \frac{1}{p}} }}} \right) = \frac{{(p - 1)q}}{p}\frac{{\left( {\sin _{p,q} x} \right)^{q - 1} \left( {1 - k^q } \right)}}{{\left( {1 - k^q \left( {\sin _{p,q} x} \right)^q } \right)^{2 - \frac{1}{p}} }},
\end{equation}
Now by using the definition and \eqref{equation2004a}, we have
\begin{eqnarray*}
 \frac{d}{{dr}}{K_{p,q} (k)}& = &\int_0^{\frac{{\pi _{p,q} }}{2}} {\frac{{r^{q-1} }}{{\left( {r'} \right)^q }}\sin _{p,q} x\frac{d}{{dx}}\left( {\frac{{ - \left( {\cos _{p,q} x} \right)^{p - 1} }}{{\left( {1 - r^q \left( {\sin _{p,q} x} \right)^q } \right)^{1 - \frac{1}{p}} }}} \right)dx}  \\
&  =& \frac{{r^{q-1} }}{{\left( {r'} \right)^q }}\int_0^{\frac{{\pi _{p,q} }}{2}} {\frac{{\left( {\cos _{p,q} x} \right)^p }}{{\left( {1 - r^q \left( {\sin _{p,q} x} \right)^q } \right)^{1 - \frac{1}{p}} }}dx}  \\
 & =& \frac{{r^{q-1} }}{{\left( {r'} \right)^q }}\int_0^{\frac{{\pi _{p,q} }}{2}} {\frac{{1 - \left( {\sin _{p,q} x} \right)^q }}{{\left( {1 - r^q \left( {\sin _{p,q} x} \right)^q } \right)^{1 - \frac{1}{p}} }}dx}  \\
 & =& \frac{{r^{ - 1} }}{{\left( {r'} \right)^q }}\int_0^{\frac{{\pi _{p,q} }}{2}} {\frac{{1 - r^q \left( {\sin _{p,q} x} \right)^q  + r^q  - 1}}{{\left( {1 - r^q \left( {\sin _{p,q} x} \right)^q } \right)^{1 - \frac{1}{p}} }}dx}  \\
&  =& \frac{1}{{r\left( {r'} \right)^q }}\left( {E_{p,q} (r) - \left( {r'} \right)^q K_{p,q} (r)} \right). \\
 \end{eqnarray*}
Similarly, we get
\begin{eqnarray*}
 \frac{d}{{dr}}{E_{p,q} (r)}&=& \int_0^{\frac{{\pi _{p,q} }}{2}} {\frac{1}{p}} \left( {1 - r^q \left( {\sin _{p,q} x} \right)^q } \right)^{\frac{1}{p} - 1} \left( { - q} \right)r^{q - 1} \left( {\sin _{p,q} x} \right)^q dx \\
&  =& \frac{q}{{pr}}\left[ {\int_0^{\frac{{\pi _{p,q} }}{2}} {\left( {1 - r^q \left( {\sin _{p,q} x} \right)^q } \right)^{\frac{1}{p}} dx - } \int_0^{\frac{{\pi _{p,q} }}{2}} {\left( {1 - r^q \left( {\sin _{p,q} x} \right)^q } \right)^{\frac{1}{p} - 1} dx} } \right] \\
  &=& \frac{q}{{pr}}\left( {E_{p,q} (r) - K_{p,q} (r)} \right). \\
 \end{eqnarray*}
With the use of formula \eqref{Kd} and \eqref{Ed}, a lengthy and trivial computation yields the proof of \eqref{Kdd} and \eqref{Edd}.$\hfill\square$

\vspace{.3cm}

For $p,q>1$ and $r\in(0,1)$, it is easy to observe that the functions $K_{p,q}(r)$ and $E_{p,q}(r)$ satisfy the following hypegeometric differential equations
$$
\frac{{d^2 }}{{dr^2 }}{K_{p,q} (r)}  - \left( {\frac{{q\left( {r'} \right)^q }}{p} + qr^q  - 2\left( {r'} \right)^q } \right)E_{p,q} (r) - \left( {2\left( {r'} \right)^{2q}  - \frac{q}{p}\left( {r'} \right)^q } \right)K_{p,q} (r)=0,
$$

$$
\frac{{d^2 }}{{dr^2 }}\left( {E_{p,q} (r)} \right) + \frac{1}{r}\left( {2 - \frac{q}{p}} \right)\frac{d}{{dr}}\left( {E_{p,q} (r)} \right) + \frac{q}{p}\frac{{r^{q - 2} }}{{1 - k^q }}E_{p,q} (r) = 0,
$$
respectively.

\vspace{.3cm}

\noindent{\bf Proof of Theorem \ref{thmturan}.}
We define 
$$f(r)=(1-r^q\sin_{p,q}(x)^q)^{1-1/p}$$
for $p,q>1,\,r\in(0,1)$ and $x\in(0,\pi_{p,q}/2)$.
Differentiating with respect to $r$, we get
$$(\log f(r))'=q\left(1-\frac{1}{p}\right)\sin_{p,q}(x)^q
\frac{r^{q-1}}{1-r^q\sin_{p,q}(x)^q}>0,$$
$$(\log f(r))''=q\left(1-\frac{1}{p}\right)\sin_{p,q}(x)^q
\frac{(q-1)r^{q-2}+r^{2(q-1)}\sin_{p,q}(x)^q}{(1-r^q\sin_{p,q}(x)^q)^2}>0.$$ 
By using the fact that the integral preserves the monotonicity and log-convexity, This implies that for $p,q>1$ and $x\in(0,\pi_{p,q}/2)$ the function 
$r\mapsto K_{p,q}(r)$ is strictly increasing and $\log$-convex on 
$(0,1)$. 

For the proof of geometrical convexity, by simple computation we get
$$f'(r)=q\left(1-\frac{1}{p}\right)\sin_{p,q}(x)^q
r^{q-1}(1-r^q\sin_{p,q}(x)^q)^{1/p-2}>0,$$
and
\begin{eqnarray*}
\left(\frac{rf'(r)}{f(r)}\right)'&=&\left(1-\frac{1}{p}\right)q\sin_{p,q}(x)^q
\left(\frac{r^q(1-r^q\sin_{p,q}(x)^q)^{1/p-2}}{1-r^q\sin_{p,q}(x)^q)^{1/p-1}}\right)'\\
&=& \left(1-\frac{1}{p}\right)q^2\sin_{p,q}(x)^q\frac{r^{q-1}}
{(1-r^q\sin_{p,q}(x)^q)^{2}}>0.
\end{eqnarray*}
Applying Lemma \ref{lemma2}, we obtain that the function 
$r\mapsto K_{p,q}(r)$
is strictly geometrically convex on $(0,1)$.

For the proof of part (2), let
$$g(r)=(1-r^q\sin_{p,q}(x)^q)^{1/p},$$
for $p,q>1,\,r\in(0,1)$ and $x\in(0,\pi_{p,q}/2)$.
A simple computation yields 
$$(\log(g(r)))'=-\frac{q}{p}\frac{r^{q-1}\sin_{p,q}(x)^q}{1-r^q\sin_{p,q}(x)^q}<0,$$
and
\begin{eqnarray*}
\left(\frac{rg'(r)}{g(r)}\right)'&=&
-\frac{q}{p}\sin_{p,q}(x)^q\frac{qr^{q-1}(1-r^q\sin_{p,q}(x)^q)-r^q(-qr^{q-1}\sin_{p,q}(x)^q)}{(1-r^q\sin_{p,q}(x)^q)^2}\\
&=& -\frac{q^2}{p}\sin_{p,q}(x)^q\frac{r^{q-1}}{(1-r^q\sin_{p,q}
(x)^q)^2}<0.
\end{eqnarray*}
Now the rest of proof follows immediately from Lemma \ref{lemma2}.
$\hfill\square$
\vspace{.3cm}

\noindent{\bf Proof of Theorem \ref{thm2806}.} For $p,q>0$, we define
$$f(p,q)=(1-t^q)^{-1/p}(1-r^qt^q)^{-1/p}$$
and 
$$g(p,q)=(1-t^q)^{-1/p}(1-r^qt^q)^{1/p}.$$
An easy computation yields that

$$\frac{\partial}{\partial q}(\log f(p,q))=\frac{t^q (\log(t)(1+r^q-2r^qt^q)+r^q\log(r)(1-t^q))}{p
   \left(1-t^q\right) \left(1-r^q t^q\right)}<0,$$

$$\frac{\partial^2}{\partial q^2}(\log f(p,q))=\frac{1}{p}\left(\frac{(\log
   (t)+\log (r))^2}{\left(1-r^qt^q\right)^2}-\frac{(\log (t)+\log (r))^2}{1-r^q t^q}+\frac{t^q \log
   ^2(t)}{\left(1-t^q\right)^2}\right)>0.$$
	This implies that for fixed $p>0$ and $r\in(0,1)$, the function
	$q\mapsto f(p,q)$ is strictly decreasing and $\log$-convex on
	$(0,\infty)$.
	
	In order to prove the monotonicity of the function $q\mapsto g(p,q)$ given in part (2), first we show that for $t\in(0,1)$ 
	and $q>0$ the function 
	$$h_q(t)=\left(\frac{t^q}{1-t^q}\right)\log\left(\frac{1}{t}\right)$$
	 has the derivative with respect to $t$
	$$(h_q(t))'=\frac{t^{q-1}(\log(1/t^q)+t^q-1)}{(1-t^q)^2},$$
	which is positive by the inequality $\log(x)>1-1/x,\,x>1$. This implies that
	the function $h_q(t)$ is strictly increasing in $t\in(0,1)$.
	Partially differentiating $ g(p,q)$ with respect $q$, we get
	$$\frac{\partial}{\partial q}(\log g(p,q))=\frac{1}{p}
	\left(h_q(rt)-h_q(t)\right),$$
	which is negative because the function $h_q(t)$ is strictly increasing in $t\in(0,1)$. Hence $g(p,q)$ is strictly decreasing 
	on $(0,\infty)$.
	
For proving the $\log$-convexity of $q\mapsto g(p,q)$, we get
	$$\frac{\partial^2}{\partial q^2}(\log g(p,q))=\frac{1}{p^2}
	(j_q(t)-j_q(rt)),$$
	where $j_q(t)=\left(\frac{t^q}{(1-t^q)^2}\right)\log(t)^2$.
	Now it is enough to prove that $j_q(t),t\in(0,1)$ is strictly increasing and positive.
	Writing $k_q(t^q)=2-2t^q+(1+t^q)\log(t^q)$, we get 
	$$j_q(t)'=\frac{t^{q-1}\log(t)k_q(t^q)}{(1-t^q)^3},$$
	which is positive, because the function 
	$k_q(t)$ is strictly increasing in $t$ and negative with 
	$k_p(1)=0$, in fact
	$$k_q(t)'=\log(t)+1/t-1>0.$$ 
	So far we have proved that the function $q\mapsto g(p,q)$ is strictly decreasing and $\log$-convex on $(0,\infty)$. 
	By using the fact that integral preserves the monotonicity and $\log$-convexity, we conclude that for fixed $p>0$ and $r\in(0,1)$ the function $q\mapsto E_{p,q}(r)$ is strictly decreasing and 
	$\log$-convex on $(0,\infty)$. This completes the proof.
$\hfill\square$

\vspace{.3cm}

\noindent{\bf Proof of Theorem \ref{yinthm}.} Let
$\alpha=r^q, \eta=q, \xi=1-\frac{1}{p}, g(x) = \left( {1 - x^q } \right)^{ - \frac{1}{p}}$. Applying formula \eqref{areq} and 

$$\frac{{\left( {\frac{1}{q}} \right)_j }}{{j!}} = ( - 1)^j 
  {- \frac{1}{q} \choose j}
 ,\quad \frac{{\left( {\frac{1}{q} + 1 - \frac{1}{p}} \right)_j }}{{j!}} = ( - 1)^j 
 {\frac{1}{p} - 1 - \frac{1}{q} \choose j },$$
we get

\begin{eqnarray*}
 b_j (\alpha ,\eta ) &=& r^{qj} \int_0^1 
{t^{jq} \left( {1 - t^q } \right)^{ - \frac{1}{p}} dt}  \\
  &=& r^{qj} \frac{1}{q}\int_0^1 {u^j \left( {1 - u} \right)^{ - \frac{1}{p}} u^{\frac{1}{q} - 1} du}  \\
  &=& \frac{{r^{qj} }}{q}B\left( {\frac{1}{q} + j,1 - \frac{1}{p}} \right) \\
  &=& \frac{{\pi _{p,q} r^{qj} }}{2}
  {- \frac{1}{q} \choose j}\displaystyle/
 {\frac{1}{p} - 1 - \frac{1}{q} \choose
 j}.
 \end{eqnarray*}
Now the claim follows easily if we apply the  formula 
$$\frac{\left( {1 - \frac{1}{p}} \right)_n }{n!} = ( - 1)^n 
 {\frac{1}{p} - 1 \choose n }.$$
This completes the proof of \eqref{Kse}. The proof of \eqref{Khat-se} and \eqref{eqse} follow similarly.
$\hfill\square$


\end{document}